\documentclass[english]{smfart}
\usepackage{amsmath, amssymb, amsthm, amscd,smfenum,xspace,mathrsfs,url,upref,verbatim}
\usepackage[utf8]{inputenc} 
\usepackage[T1]{fontenc}
\usepackage{dsfont}
\usepackage{mathabx}
\usepackage{stmaryrd}
\usepackage[all,cmtip]{xy}
\usepackage{xfrac} 

\usepackage{ textcomp }

\let\tilde=\widetilde
\numberwithin{equation}{subsection}

\newtheorem{theorem}{Theorem} 
\newtheorem{proposition}[equation]{Proposition}

\theoremstyle{remark}

\DeclareMathOperator{\Tor}{Tor}

\DeclareMathOperator{\Supp}{Supp}

\def\cartesien{\ar@{}[rd]|{\square}}

\DeclareMathOperator{\Irr}{Irr}

\DeclareMathOperator{\hol}{hol}

\author[J.-B.~Teyssier]{Jean-Baptiste Teyssier}
\curraddr{Freie Universität Berlin, Mathematisches Institut, Arnimallee 3, 14195 Berlin, Germany}
\email{teyssier@zedat.fu-berlin.de}

\title{
Tensor product and irregularity for holonomic $\mathcal{D}$-modules}
\usepackage[french,english]{babel}
\begin{document}
\maketitle
\section*{Introduction}
Let $X$ be a complex variety and let $D^{b}_{\hol}(\mathcal{D}_X)$ be the derived category of complexes of $\mathcal{D}_X$-modules with bounded holonomic cohomology. It is known \cite[6.2-4]{Mehbsmf} that for a regular complex\footnote{that is, a complex with regular cohomology modules.} $\mathcal{M}\in D^{b}_{\hol}(\mathcal{D}_X)$, the derived tensor product $\mathcal{M}\otimes^{\mathds{L}}_{\mathcal{O}_X}\mathcal{M}$ is regular. The goal of this note is to prove the following
\begin{theorem}\label{maintheorem}
Let $\mathcal{M}\in D^{b}_{\hol}(\mathcal{D}_X)$ and suppose that $\mathcal{M}\otimes^{\mathds{L}}_{\mathcal{O}_X}\mathcal{M}$ is regular. Then $\mathcal{M}$ is regular. 
\end{theorem}
The technique used in this text is similar to that used in \cite{TeyReg}, and proceed by recursion on the dimension of $X$. The main tool is a sheaf-theoretic measure of irregularity \cite{Mehbgro}. \\ \indent

I thank the anonymous referee for a careful reading and constructive remarks on this manuscript. This work has been achieved with the support of Freie Universität/Hebrew University of Jerusalem joint post-doctoral program and ERC 226257 program.
\subsection{}For every morphism $f:Y\longrightarrow X$ with $X$ and $Y$ complex varieties, we denote by $f^{+}:D^{b}_{\hol}(\mathcal{D}_X)\longrightarrow D^{b}_{\hol}(\mathcal{D}_Y)$ and $f_{+}:D^{b}_{\hol}(\mathcal{D}_Y)\longrightarrow D^{b}_{\hol}(\mathcal{D}_X)$ the inverse image and direct image functors for $\mathcal{D}$-modules. We define $f^{\dag}:=f^{+}[\dim Y-\dim X]$.
\subsection{}\label{defirr}
If $Z$ is a closed analytic subspace of $X$, we denote by $\Irr^{\ast}_{Z}(\mathcal{M})$ the irregularity sheaf of $\mathcal{M}$ along $Z$  \cite{Mehbgro}.
\section{The proof}

\subsection{The 1-dimensional case}
We suppose that $X$ is a neighbourhood of the origin  
$0\in \mathds{C}$ and we prove the following
\begin{proposition}\label{propdim1}
Let $\mathcal{M}\in D^{b}_{\hol}(\mathcal{D}_X)$ so that $\mathcal{H}^{k}\mathcal{M}$ is a smooth connexion away from 0 for every $k\in \mathds{N}$.  If
$\mathcal{M}\otimes^{\mathds{L}}_{\mathcal{O}_X}\mathcal{M}$
is regular, then $\mathcal{M}$ is regular.
\end{proposition} 
 The complex $$(\mathcal{M}\otimes^{\mathds{L}}_{\mathcal{O}_X}\mathcal{M})(\ast 0)\simeq 
\mathcal{M}(\ast 0)\otimes_{\mathcal{O}_X}\mathcal{M}(\ast 0)
$$
is regular. Since we are in dimension one, the regularity of $\mathcal{H}^{k}\mathcal{M}$ is equivalent to the regularity of   $\mathcal{H}^{k}\mathcal{M}(\ast 0)$. Thus, we can suppose that $\mathcal{M}$ is localized at $0$. In particular, the $\mathcal{H}^{k}\mathcal{M}$  are  flat $\mathcal{O}_X$-modules, so the only possibly non zero terms in the Künneth spectral sequence
\begin{equation}\label{ss1}
E_{2}^{pq}=\bigoplus_{i+j=q}\Tor^{p}_{\mathcal{O}_X}(\mathcal{H}^{i}\mathcal{M},\mathcal{H}^{j}\mathcal{M})  \Longrightarrow \mathcal{H}^{p+q}(\mathcal{M}\otimes_{\mathcal{O}_X}\mathcal{M})
\end{equation}
sit on the line $p=0$. Hence, \eqref{ss1} degenerates at page 2 and induces a canonical identification
$$
\mathcal{H}^{k}(\mathcal{M}\otimes_{\mathcal{O}_X}\mathcal{M})\simeq \bigoplus_{i+j=k}(\mathcal{H}^{i}\mathcal{M}\otimes_{\mathcal{O}_X}\mathcal{H}^{j}\mathcal{M})
$$
 for every $k$. In particular, the module $\mathcal{H}^{i}\mathcal{M}\otimes_{\mathcal{O}_X}\mathcal{H}^{i}\mathcal{M}$ is regular for every $i$. Hence, one can suppose that $\mathcal{M}$ is a germ of meromorphic connexions at 0. By  looking formally at $0$, one can further suppose that $\mathcal{M}$ is a $\mathds{C}((x))$-differential module. In this case, the regularity of $\mathcal{M}$ is a direct consequence of the Levelt-Turrittin decomposition theorem \cite{SVDP}.
\subsection{Proof of theorem \ref{maintheorem} in higher dimension}
We proceed by recursion on the dimension of $X$ and suppose that $\dim X>1$. For every point $x\in X$ taken away from a discrete set of points $S\subset X$, one can find a smooth hypersurface $i:Z\longrightarrow X$ passing through $x$ which is  non characteristic  for $\mathcal{M}$. Since regularity is preserved by inverse image, the complex
$$
i^{+}\mathcal{M}\otimes^{\mathds{L}}_{\mathcal{O}_X}i^{+}\mathcal{M}
$$
is regular. By recursion hypothesis, we deduce that $i^{+}\mathcal{M}$ is regular.
From \cite[3.3.2]{TeyReg}, we obtain
$$
\Irr^{\ast}_{x}(\mathcal{M})\simeq \Irr^{\ast}_{x}(i^{+}\mathcal{M})\simeq 0
$$
Since regularity can be punctually tested \cite[6.2-6]{Mehbsmf}, we deduce that $\mathcal{M}$ is regular away from $S$.  In what follows, one can thus suppose that $X$ is a neighbourhood of the origin $0\in \mathds{C}^n$
and that  $\mathcal{M}$ is regular away from 0.
\\ \indent
Let us suppose that 0 is contained in an irreducible component of $\Supp \mathcal{M}$ of dimension $>1$. Let $Z$ be an hypersurface containing 0 and satisfying the conditions 
\begin{enumerate}
\item $Z\cap \Supp\mathcal{M}$ has codimension $1$ in $\Supp\mathcal{M}$.
\item The modules  $\mathcal{H}^{k}\mathcal{M}$ are smooth\footnote{That is, $\Supp(\mathcal{H}^{k}\mathcal{M})$ is smooth away from $Z$ and the characteristic variety of $\mathcal{H}^{k}\mathcal{M}$ away from $Z$ is the conormal bundle of  $\Supp(\mathcal{H}^{k}\mathcal{M})$ in $X$.} away from $Z$.
\item 
$\dim \Supp R\Gamma_{[Z]}\mathcal{M}<\dim \Supp\mathcal{M}$.
\end{enumerate}
Such an hypersurface always exists by \cite[6.1-4]{Mehbsmf}. According to the fundamental criterion for regularity \cite[4.3-17]{Mehbsmf}, the complex   $\mathcal{M}(\ast Z)$ is regular. From the local cohomology triangle 
$$
\xymatrix{
R\Gamma_{[Z]}\mathcal{M}\ar[r] &   \mathcal{M} \ar[r] &   \mathcal{M}(\ast Z)  \ar[r]^-{+1}&
}
$$
we deduce that one is left to prove that  $R\Gamma_{[Z]}\mathcal{M}$ is regular. There is a canonical isomorphism
\begin{equation}\label{isocan}
R\Gamma_{[Z]}\mathcal{M}\otimes^{\mathds{L}}_{\mathcal{O}_X}
R\Gamma_{[Z]}\mathcal{M} \simeq R\Gamma_{[Z]}(\mathcal{M}\otimes^{\mathds{L}}_{\mathcal{O}_X}\mathcal{M} )
\end{equation}
Since $R\Gamma_{[Z]}$ preserves regularity, the left hand side of \eqref{isocan} is regular. So one is left to prove theorem \ref{maintheorem} for $R\Gamma_{[Z]}\mathcal{M}$, with $\dim \Supp R\Gamma_{[Z]}\mathcal{M}<\dim \Supp \mathcal{M}$. By iterating this procedure if necessary,  one can suppose that the components of  $\Supp \mathcal{M}$ containing 0 are curves. We note $C:=\Supp \mathcal{M}$. At the cost of restricting the situation to a small enough neighbourhood of 0, one can suppose that $C$ is smooth away from $0$. Let $p:\tilde{C}\longrightarrow X$ be the composite of normalization map for $C$ and the canonical inclusion $C\longrightarrow X$. By Kashiwara theorem \cite[1.6.1]{HTT}, the canonical adjunction \cite[7.1]{MS}
\begin{equation}\label{isocan2}
p_{+} p^{\dag}\mathcal{M}\longrightarrow 
\mathcal{M}
\end{equation}
is an isomorphism away from 0. So the cone of \eqref{isocan2} is supported at $0$. Hence, it is regular. One is then left to show that $p_{+} p^{\dag}\mathcal{M}$ is regular. Since regularity is preserved by proper direct image, we are left to prove that $p^{\dag}\mathcal{M}$ is regular. There is a canonical isomorphism
\begin{equation}\label{isocan3}
p^{\dag}\mathcal{M} \otimes^{\mathds{L}}_{\mathcal{O}_{\tilde{C}}}p^{\dag}\mathcal{M} 
\simeq p^{\dag}(\mathcal{M} \otimes^{\mathds{L}}_{\mathcal{O}_X}\mathcal{M} )
\end{equation}
So the left hand side of \eqref{isocan3} is regular and one can apply \ref{propdim1}, which concludes the proof of theorem \ref{maintheorem}.
\bibliographystyle{amsalpha}
\bibliography{KL}

\end{document}